\documentclass[a4paper,12pt]{article}

\usepackage[cp1251]{inputenc}
\usepackage[russian]{babel}
\usepackage{enumerate}

\newtheorem{theorem}{\it Теорема}
\newtheorem{lemma}{\it Лемма}

\title{Конечность продолжительности карточной игры <<Разори моего соседа>>}

\author{Е.Лакштанов\thanks{Department of Mathematics, Aveiro University, Aveiro 3810, Portugal.  This work was supported by {\it FEDER} founds through {\it COMPETE}--Operational Programme Factors of Competitiveness (``Programa Operacional Factores de Competitividade'') and by Portuguese founds through the {\it Center for Research and Development in Mathematics and Applications} (University of Aveiro) and the Portuguese Foundation for Science and Technology (``FCT--Fund\c{c}\~{a}o para a Ci\^{e}ncia e a Tecnologia''), within project PEst-C/MAT/UI4106/2011 with COMPETE number FCOMP-01-0124-FEDER-022690 (lakshtanov@rambler.ru).} \and
 Алена Алексенко$^*$}

\begin{document}
\maketitle

\begin{abstract}Для карточных игр вида <<Разори Моего Соседа>> доказана конечность математического ожидания продолжительности игры при условии, что начинающий розыгрыш игрок, определяется случайно, а взятка перемешивается при добавлении в колоду.
Результат верен для модификаций правил игры общего вида.
\end{abstract}



В \cite{lr} мы доказали конечность математического ожидания количества ходов в карточной игре  <<Пьяница>>, для различных ее модификаций. Напомним, как выглядит наиболее простая версия этой  игры: колода состоит из карт одной масти и в начале игры делится на две равные части. Каждый игрок хранит свои карты в виде колоды, в которой все карты расположены рубашками вверх. Игроки снимают по одной верхней карте, и обе возвращаются в низ колоды игрока, карта которого оказалось старше. О разновидностях этой игры можно прочитать в \cite{wiki2}. В основном, модификации этой игры относятся к количеству карт в колоде, порядку старшинства (в русской версии самая младшая карта старше чем <<туз>>), а также к порядку действий, если в ходе розыгрыша игроки предъявят карты одинакового достоинства. Игру в <<Пьяницу>> естественным образом можно рассматривать как марковскую цепь, поскольку в колоду возвращается набор карт, порядок которого  не определяется правилами.

Давшая название одноименному экономическому термину \cite{wiki3} игра <<Beggar my Neighbour>> (Разори моего соседа) упоминалась еще в романе Чарльза Диккенса в 1860 году \cite{dick}. Ее правила являются детерминированными: в начале игры каждый из двух игроков получает  половину колоды карт, которые хранятся как отдельные колоды <<рубашками>> вверх и карты не перемешиваются в течение игры.  В ходе розыгрыша игрок может брать только верхнюю карту своей колоды. Игроки по очереди кладут по карте на стол до тех пор пока не будет выложена <<картинка>>, то есть <<валет>>, <<дама>>, <<король>> или <<туз>>: в  этом случае, другой игрок платит <<штраф>> выкладывая одну за одной определенное количество своих карт: одну карту в случае <<валета>>, две в случае <<дамы>>, три в случае <<короля>> и четрые в случае <<туза>>. Если в процессе выплаты штрафа на стол выкладывается <<картинка>>, то процесс выплаты останавливается, и уже другой игрок начинает платить штраф. Если штраф выплачен полностью, то игрок (которому платили штраф) забирает всю взятку и размещает снизу  своей колоды. Как видно из описания, существенно важно, кто начинал розыгрыш, то есть кто положил первую карту. {\it Игрок взявший взятку начинает следующий розыгрыш!} Игрок оставшийся без карт, проигрывает. Если игрок не может сделать ход предусмотренный правилами, то есть процедура описанная в правилах розыгрыша не может быть закончена, то этот игрок проигрывает.

Можно рассмотреть другие игры, в которых изменены правила определяющие процесс розыгрыша: можно изменить количество карт требуемых для выплаты штрафа, также как и сам набор  карт, появление которых требует штрафа. Игры, получающиеся при помощи указанной процедуры будем называть играми вида <<разори моего соседа>>. Игру с оригинальными правилами будем называть, {\it классической} игрой <<разори моего соседа>>.

Естественен вопрос о наличии циклов в этой игре, он был включен в <<анти-Гильбертовый>> список Джона Конвэя \cite{wiki},\cite{guy}. Вопрос был открытым более  47 лет пока  Marc Paulhus не сумел  обнаружить циклы \cite{con} для стандартной колоды в 52 карты. Джон Конвэй пишет так об этом: "About 1 in 150.000 games goes on forever. Just as surely, however, the total number of times this game has been played by the World's children must be significantly larger than 150.000, so many of them will have been theoretically non-terminating ones. We imagine, though, that in practice most of them, actually did terminate because someone made a mistake..."

Нам удалось доказать конечность количества ходов требуемых для окончания игры, если взятка перемешивается при возвращении в колоду, и, что самое важное, если игрок, начинающий новый розыгрыш  определяется испытанием Бернулли с $p\in(0,1)$. Видимо, так и происходит на самом деле, что иногда игроки начинают розыгрыш не соблюдая инструкции, так что наша модель претендует на описание реального процесса.
Подобная общность является <<недостатком>> нашего метода, с другой стороны, его преимуществом является, простота рассуждений и применимость для правил общего вида.

\section{Общий вид правил}
Пусть $W$ это конечное множество карт. Множеством состояний в игре является произвольное разбиение $W$ на два упорядоченных подмножества:
$$
\mathcal H = \{ (L,R) ~ : ~ L \sqcup R=W\}
$$
Переход от одного состояния к другому определяется некоторым процессом, называемым розыгрышем. Игроки берут карты сверху своей колоды и выкладывают на стол. Правила игры определяют по упорядоченному набору карт выложенному на стол, что будет происходить дальше: либо игра продолжается и кто-то кладет следующую карту на стол, либо розыгрыш закончен и один из игроков забирает всю взятку себе, размещая снизу своей колоды.

В этой работе мы считаем, что каждый раз, игрок начинающий розыгрыш определяется случайно, как  результат испытания Бернулли с вероятностью $p \in (0,1)$.
Итак правила игры определяются некоторой функцией $F$ заданной на множестве $S$ всех упорядоченных подмножеств $W$:
$$
F ~ : ~ S \rightarrow \{\mbox{<<Закончить>>},\mbox{<<Продолжить>>}\} \times \{1,2\}
$$

Пусть $S_1$ есть множество состоящее из одной верхней карты игрока номер которого есть результат испытания бернулли. Если на столе были выложен упорядоченный набор $S_{n} \subset S$, то значение $F(S_n)=(\mbox{<<Продолжить>>} ~ \!, ~ i), ~ i=1,2$ означает, что $i-\!$ый игрок выкладывает на стол верхнюю карту $t \in W$
$$
S_{n+1}=\{t\} \amalg S_n,
$$
 а если $F(S_n)=(\mbox{<<Закончить>>} ~ \!, ~i), ~ i=1,2$ то  $i-\!$ый игрок забирает взятку, добавляя перетасованные карты в конец своей колоды. Символ $\amalg$ применяемый к двум упорядоченным множеством, создает новое упорядоченное множество, состоящее из объединения двух исходных, сохраняя исходный порядок множеств, причем вначале перечисляются элементы первого аргумента $\amalg$. После окончания розыгрыша проводится новое испытание бернулли и выигравший его начинает новый розыгыш.  Нам потребуется также наложить условие {\it невырожденности для игр <<разори моего соседа>>}, неформальный смысл которого состоит в следующем. Рассмотрим произвольную раздачу карт $(L,R)$ и проследим за двумя розыгрышами для каждого из двух возможных начинающих розыгрыш игроков. Условие состоит в том, что пары (<<номер игрока берущего взятку>> И <<упорядоченный набор карт при окончании розыгрыша>>) должны быть различными. Сформулируем это условие на формальном языке.

Обозначим через $\mathcal G$ граф, вершинами которого являются произвольные раздачи карт $(L,R) \in \mathcal H$, а набором ребр являются упорядоченные пары вершин, соответствующие переходу согласно процедуре, которую мы сейчас опишем:

Намомним, что мы обозначили через $S$ множество всех упорядоченных наборов $W$. Определим операторы $M_i$ действующие на произвольных тройках $(L,R,S)$ таких что
 $$
 L \sqcup R \sqcup S =W
 $$
 следующим образом
 $$
 M_1 (t \amalg L, R,S)=(L,R,t \amalg S), \quad  M_2 ( L, t \amalg R,S)=(L,R,t \amalg S),
 $$
 то есть оператор $M_i$ соответствует тому, что $i-$ый игрок положил в ходе розыгрыша свою верхнюю карту на стол, причем сверху карт уже выложенных в течение этого розыгрыша.

Поскольку каждый из игроков может начать розыгрыш, каждому состоянию $(L_0,R_0)$ сопоставим две тройки
$$
(L^i_1,R^i_1,S^i_1) := M_i(L_0,R_0,\emptyset), \quad i=1,2.
$$
Индекс $i$ соответствует номеру игрока, начинающему розыгрыш.
Далее если $F(S_n^i)=(\mbox{<<продолжить>>} ~ \!, ~m), ~ m=1,2$ то
\begin{equation}\label{actionM}
(L^i_{n+1},R^i_{n+1},S^i_{n+1}) := M_{m}(L^i_n,R^i_n,S^i_n),
\end{equation}
то есть розыгрыш продолжается и следующую карту кладет игрок с индексом $m$.
Если же $F(S_n^i)=(\mbox{<<Закончить>>} ~ \!, ~m), ~ m=1,2$ то в случае $m=1$ из вершины $(L_0,R_0) \in \mathcal G$ выходят ребра во все вершины вида  $(L^i_n \amalg \widehat{S}^i_n,R^i_n)$, где $\widehat{S}^i_n$ есть произвольная перестановка $S^i_n$. Среди всех возможных потомков, мы выделим те, которые позволяют однозначно определить предка: обозначим $(L_0,R_0)^i:=(L^i_n \amalg S^i_n,R^i_n)$.

По аналогии, если $m=2$, то $(L_0,R_0)$ переходит в одно из состояний $(L^i_n,R^i_n \amalg \widehat{S}^i_n)$. Также обозначим $(L_0,R_0)^i:=(L^i_n,R^i_n \amalg {S}^i_n)$.

Вершины, соответствующие концу игры, не имеют выходящих ребер.  Такие вершины называются абсорбирующими.

{\bf Условие невырожденности для игр <<разори моего соседа>>.} Мы предполагаем, что набор карт $W$ и правила $F$ таковы, что если для некоторой раздачи карт $(L_0,R_0)$ верно, что
$$
(L_0,R_0)^1 =(L_0,R_0)^2,
$$
то состояние
$(L_0,R_0)^1 =(L_0,R_0)^2$ соответстует концу игры.

Простыми словами, результат любого розыгрыша всегда зависит от того, кто из двух игроков положил свою карту первым.
\\
{\it Замечание 1.} Легко видеть, что игра описанная в начале <<Разори моего Соседа>> удовлетворяет условию невырожденности. Действительно, допустим, что для  какого-либо состояния
$(L_0,R_0)$ выполняется $(L_0,R_0)^1 = (L_0,R_0)^2$. В частности, это означает что наборы карт выложенные в каждом случае, совпадают. То есть либо карты заканчиваются (конец игры) либо один из игроков берет взятку, но, в силу симметрии, номер игрока зависит от того, кто начинал розыгрыш, следовательно $(L_0,R_0)^1 \neq (L_0,R_0)^2$ потому, что в одном случае набор $S^1_n=S^2_n$ присоединяется к левому множеству, а в другом случае, к правому.

{\it Замечание 2.} Игра <<пьяница>>, очевидно, не удовлетворяет условию невырожденности, поскольку игроки могут показывать свои карты одновременно.

\begin{theorem}\label{main}
Если набор карт $W$ и правила $F$ таковы, что выполнено условие невырожденности игры, то математическое ожидание количества ходов требуемых для конца игры конечно.
\end{theorem}

В силу замечания 1 имеем \\
{\it Следствие 1.}
Аналог игры <<Разори моего соседа>>, в котором игрок, начинающий каждый розыгрыш определяется испытанием Бернулли с вероятностью $p\in (0,1)$, а взятка перемешивается при возвращении в колоду, имеет конечное математическое ожидание количества ходов требуемых для окончания игры.

{\bf Замечание. }Заметим, что переход от $(L_0,R_0)$ к $(L_0,R_0)^i$ не соответствует переходу в классической игре (без перемешивания взятки), но соответствовал бы, если играть  так, что колода  карт, у каждого игрока, лежит рубашками вниз,  то есть игроки видят значения верхних карт.  При таком способе игры, результат о конечности будет верен, без условия на то, что карты должны перемешиваться при возвращении в колоду.

\section{Доказательство}
В \cite{lr} мы подробно обсуждали, что математическое ожидание количества ходов до абсорбции конечно тогда и только тогда, когда граф цепи является абсорбирующим, то есть от каждой вершниы есть хотя бы один путь ведет до абсорбирующей вершины. Для того, чтобы показать абсорбируемость графа $\mathcal G$ мы воспользуемся достаточным условием абсорбируемости, полученным в \cite{lr}.
\begin{lemma}\cite{lr}\label{lemmaold}
Если сущесвтует подграф $G_0 \subset G$ с тем же множеством вершин, что и $G$ обладающий следующими свойствами: Для каждой неабсорбирующей вершины $v \in G_0$ верно
\begin{itemize}
\item Выходящая степень $deg^{out} v$ равна 2.
\item Входящая степень $deg^{in} v$  не превосходит 2.
\item Если $deg^{in} v =2$ то у нее есть предок $w \in G_0$ такой, что $deg^{in} w <2 $.
\end{itemize}
(степень берется в $G_0$). Тогда граф $G$ является абсорбирующим.
\end{lemma}
Доказательство этой леммы несложно, поэтому повторим его здесь. Будем называть вершину {\it достигающей} если среди ее потомков есть абсорбирующая вершина, {\it блуждающей} в обратном случае. Очевидно, что потомок блуждающей вершины является блуждающей вершиной, а предок достигающей является достигающей вершиной. Обозначим через $\mathcal W$ подмножество всех достигающих вершни графа $\mathcal G_0$. Суммарное количество ребер выходящих из вершин $\mathcal W$ есть в точности удвоенное количество вершин в $\mathcal W$, а поскольку ребра вышедшие из вершин $\mathcal W$ возвращаются в это же множество, и по условие леммы $deg^{in} v \leq 2$  мы получаем, что входящая степень блуждающей вершины есть 2, и вершины $\mathcal W$ образуют изолированный подграф, то есть в него нельзя попасть извне. Поэтому, при выполнении предположения 3, множество $\mathcal W$ является пустым. Лемма доказана.

Пусть $\mathcal G_0$ есть подграф $\mathcal G$, с тем же набором вершин, но с меньшим количеством ребер. Для каждой вершиниы $(L_0,R_0)$ оставим только два ребра соответствующие $(L_0,R_0)^1,(L_0,R_0)^2$. Покажем, что у каждой вершины $\mathcal G_0$ не более двух прямых предков. Рассмотрим карты снизу колоды первого игрока. По построению, последняя карта в колоде, являлась первой в розыгрыше, приведшем к взятию взятки первого игрока. Предпоследняя карта являлась второй в розыгрыше и так далее. Итак, по набору карт снизу можно однозначно определить как происходил розыгрыш, и однозначно построить прямого предка соответствующего тому, что взятку забирал первый игрок.  Также строится прямой предок, соответствующий тому, что взятку на предыдущем ходу забирал второй игрок. Будем использовать терминологию <<левый>> прямой предок и <<правый>> прямой предок.  Если в процедуре восстановления предка, карты в колоде заканчивается раньше, то это состояние не является по построению, прямым потомком какой-либо вершины, и входящая степень меньше двух. Если же и левый и правый предок существуют, то они различны в силу условия невырожденности, и в этом случае входящая степень равна двум. Заметим, что переход к левому предку соответствует тому, что количество карт у  первого игрока уменьшается, поскольку часть карт передается второму игроку. Аналогично переход к правому предку соответствует тому, что количество карт у второго игрока уменьшается.

Для удовлетворения условиям леммы осталось показать, что у любой вершины, с выходящей степенью два, есть предок с выходящей степенью меньшей чем два. Рассмотрим вершину $v=(L_0,R_0), ~ deg^{out}v=2$. Переход к левому предку влечет за собой уменьшение количества карт у первого игрока, что не может продолжать бесконечное количество раз, следовательно есть предок у которого не более одного прямого предка, так как по построению, вершины соответствующие концу игры не имеют выходящих вершин.

Итак, мы показали, что граф $\mathcal G_0 \subset \mathcal G$ удовлетворяет условиям леммы \ref{lemmaold}, а следовательно граф $\mathcal G$ является абсорбирующим. Теорема 1 доказана.

\newpage

Алена Ильинична Алексенко, университет Авейру, Португалия, a40861@ua.pt
\\
\\
Автор ответственный за переписку: Евгений Леонидович Лакштанов, университет Авейру, Португалия, lakshtanov@ua.pt, телефон: +351 234 18 66 00
почтовый адрес: University of Aveiro, department of mathematics, 3810-193, Aveiro, Portugal

\end{document}